\documentclass[10pt,a4paper]{article}
\usepackage[all]{xy}
\usepackage{theorem}
\usepackage{amssymb}
\usepackage{amsmath}
\usepackage{latexsym}

\newcommand\defi[2][]{\emph{#2}}

\newtheorem{The}{Theorem}
\newtheorem{Lem}[The]{Lemma}\newtheorem{Pro}[The]{Proposition}
\newcommand{\bevis}{\textit{Proof: }}\newcommand{\qed}{\hfill$\Box$}

\newcommand{\ned}{\vspace{\theorempostskipamount}\noindent}

\newcommand{\ra}{\rightarrow}

\newcommand{\wt}{\widetilde}
\newcommand{\ol}{\overline}\newcommand{\ul}{\underline}

\begin{document}
\title{Equations for some very special \\ 
   Enriques surfaces in characteristic two}
\author{\textsc{Pelle Salomonsson}}\date{}
\maketitle

\abstract{\noindent We give an explicit construction and 
clasification of some very special sort of Enriques surfaces 
in characteristic two. This proves the 
existence of some of the surfaces 
that were called  ``extra-special'' by Cossec and Dolgachev in 
their book on the subject, and whose existence was left open by them. 
These surfaces coincide with a class of surfaces that are introduced 
in recent work by Ekedahl and Shepherd--Barron. }

\section{Introduction and statement of results}
A minimal and smooth algebraic surface~$X$ is called an 
\defi{Enriques} surface if 
$\chi(X,\mathcal{O}_{X})=1$ and 
$\omega^{\otimes2}_{X}\cong\mathcal{O}_{X}$, 
where $\omega_{X}$ is the canonical line bundle. We shall study them 
in characteristic two here, that is, the ground field has 
characteristic two. Enriques surfaces depend on ten moduli and 
typically they  do not contain 
any smooth rational curves. So it is natural to study those 
that have many such curves, and more precisely those surfaces having 
some specified \emph{configuration} of rational curves, intersecting 
each other according to the given intersection graph. In the 
book~\cite{CD} by Cossec and Dolgachev it is suggested that the most 
special surfaces should be those that have a configuration of smooth 
rational curves intersecting each other according to an extended 
Dynkin graph plus an additional rational curve~$R$ intersecting  
that configuration in a suitable way. 

We shall be interested in 
the extended Dynkin graphs~$\wt{E}_{6}$,~$\wt{E}_{7}$ and 
$\wt{E}_{8}$ since these are the most special in some sense and hence 
the most interesting. We shall call such a configuration an 
\defi{$(\wt{E}_{i}+R)$-configuration}, and the 
Enriques surface containing it an \defi{$(\wt{E}_{i}+R)$-surface}.  
If the curve~$R$ is omitted, then we speak about 
\defi{$\wt{E}_{i}$-configurations}. 
The components of the $\wt{E}_{i}$-configuration can be given 
multiplicities in such a manner that this non-reduced curve appears 
a half-fiber in a fibration of the surface over the projective line. 
The curve~$R$ then intersects the $\wt{E}_{i}$-configuration 
in one of its components occurring with multiplicity one. 

We must also recall the fact that in characteristic two there are
three types of Enriques surfaces according to the nature of the
numerically trivial part of their Picard group schemes, which in all
cases are finite group schemes of order two. 
We have the three cases~$\pmb{\mu}_{2}$,~$\mathbf{Z}/2$ 
and~$\pmb{\alpha}_{2}$. Only the latter two are of interest here 
since the first type cannot
have $(\wt{E}_{i}+R)$-configurations on them (as we will se). 
What sets the types apart is among other 
things that the canonical class is 
\emph{trivial} in the two cases with connected Picard schemes. More 
concretely, that means that when representing the surface as a curve 
fibration with genus-one fibers, these two types will have a 
\emph{single} non-reduced fiber of multiplicity two (a ``wild'' double 
fiber), whereas the $\mathbf{Z}/2$ case behaves just like the 
characteristic zero case, having two double fibers along which the 
canonical class is concentrated.

The overall logical structure of this moduli problem is that we have a 
set~$M$ of polynomial equations defining algebraic surfaces and we 
have another set~$M'$ of certain moduli data involving Enriques 
surfaces, and we must show $i)$~that the natural 
mapping from~$M$ to the set 
of isomorphism classes of algebraic surfaces factors through a mapping 
$M{\ra}M'$, and $ii)$~that this mapping is surjective, 
and $iii)$~that it is injective. The last point will be proved only 
in the~$\mathbf{Z}/2$ case since it is false in the other case with 
our choice of equations. The first point is easiest and will be 
considered first. 

The following theorem deals with the first of the three problems: that 
the equations define Enriques surfaces and not some other sort of 
algebraic surfaces. The equations are multihomogeneous in 
the variables $z,x,y,s,t$ such
that writing~$d_{x}$ for the degree of~$x$ et cetera we have for each
monomial $d_{s}+d_{t}+2d_{z}=4$ and
$2d_{t}+d_{x}+d_{y}+5d_{z}=10$. The variables $x,y$ are coordinates on
the projective line and $s$ and $t$ are supplementary coordinates on
the ruled surface 
$\mathbf{P}(\mathcal O_{\mathbf{P}^1}(2)s\bigoplus
\mathcal O_{\mathbf{P}^1}t)$ and
the form of the equations will make our surface a double cover of that
ruled surface. The rest of the appearing letters are not 
variables but parameters, that is, coefficients. 
The coefficient in front of $x^{i}y^js^kt^{\ell}$ is usually  
denoted~$a_{ij}$ and the one in front of $zx^{i}y^js^kt^{\ell}$ is 
usually written~$b_{ij}$. In some cases there are dependencies among 
the parameters, and we therefore have some parameters that appear as 
coefficients at several places. These are denoted~$v$ or~$w$. 

\begin{The} \label{firsttheorem}
    The following equations cut out surfaces that upon 
resolving singularities are Enriques surfaces of types 
$\mathbf{Z}/2$ or~$\pmb{\alpha}_{2}$ with an 
$(\wt{E}_{i}+R)$-configuration. For 
an explanation of the appearing symbols, see the preceding 
paragraph. 
\begin{gather*}
\ul{\text{$\mathbf{Z}/2$ type with $\wt{E}_6+R$ ($v\neq0$):}}  \\
z^2+z(b_{32}x^3y^2s^2+v^2x^2yst)+  \\
(y^4+x^4)x^3y^3s^4+
 (v^3xy+a_{53}x^2)x^3y^3s^3t+v^2x^3y^3s^2t^2+xyt^4=0, \\ 
\ul{\text{$\mathbf{Z}/2$ type with $\wt{E}_7+R$ ($b_{32}\neq0$):}} \\
z^2+zb_{32}x^3y^2s^2+ 
(y^4+x^4)x^3y^3s^4+a_{53}x^5y^3s^3t+xyt^4=0, \\ 
\ul{\text{$\mathbf{Z}/2$ type with $\wt{E}_8+R$ ($a_{53}\neq0$):}}\\
z^2+ (y^4+x^4)x^3y^3s^4+a_{53}x^5y^3s^3t+xyt^4=0, \\
\ul{\text{$\pmb{\alpha}_{2}$ type with $\wt{E}_6+R$ ($w\neq0$):}} \\
z^2 + z(wxy+b_{50}x^2)x^3s^2 + \\ x^3y^7s^4 + 
  (wy^3+a_{80}x^3)x^5s^3t+wx^4st^3 + xyt^4=0 \\
\ul{\text{$\pmb{\alpha}_{2}$ type with $\wt{E}_7+R$ ($b_{50}\neq0$):}} \\
z^2 + zb_{50}x^5s^2 +x^3y^7s^4 + a_{80}x^8s^3t+ xyt^4=0 \\
\ul{\text{$\pmb{\alpha}_{2}$ type with $\wt{E}_8+R$ ($a_{80}\neq0$):}}\\
z^2  +x^3y^7s^4 + a_{80}x^8s^3t+ xyt^4=0 \,.
\end{gather*}
\end{The}

\begin{The} \label{secondtheorem}
    $i)$ An Enriques surface of type~$\mathbf{Z}/2$ 
    or~$\pmb{\alpha}_{2}$ with an $(\wt{E}_{i}+R)$-configuration 
    can be written on the form stated in theorem~\ref{firsttheorem}. 
    
    $ii)$ In the~$\mathbf{Z}/2$ case, the surface can be 
    \emph{uniquely} written in that form. 
    
    $iii)$ In the~$\pmb{\alpha}_{2}$ case, the surface can be 
    written in that form in such a manner that the 
    parameters~$w$,~$b_{50}$ or~$a_{80}$ are equal to~1 in the 
    $\wt{E}_6$, $\wt{E}_7$ and $\wt{E}_8$ cases, respectively.   
\end{The}

\ned
In~\cite{ESB} Ekedahl and Shepherd--Barron introduce the notion of 
``exceptional Enriques surfaces''. The definition will not be 
reproduced here but let us mention that an
exceptional $\mathbf{Z}/2$-type Enriques surface always has a
non-trivial vector field. These 
authors show that an Enriques surface in characteristic two is 
exceptional if and only if it has an $(\wt{E}_n+R)$-configuration 
for $n=6,7$ or~$8$. So the above two theorems give a 
classification of such surfaces. 

Recall that an $\wt{A}_{1}$-configuration 
($\wt{A}_{1}^{*}$-configuration) is two 
smooth rational curves 
intersecting transversally in two points 
(in a single point with multiplicity two). 
On page 186 in~\cite{CD} the authors introduce the notions of 
$\wt{D}_{8}$-special, $\wt{E}_{8}$-special and 
$(\wt{A}_{1}^{}+\wt{E}_{7})$-special Enriques 
surfaces. An $\wt{E}_{8}$-special surface is 
the same as what we have called an 
$(\wt{E}_{8}+R)$-surface and $\wt{D}_{8}$-special surfaces are 
defined analogously. An $(\wt{A}_{1}^{}+\wt{E}_{7})$-special 
($(\wt{A}_{1}^{*}+\wt{E}_{7})$-special) surface 
has an $(\wt{E}_{7}+R)$-configuration and an 
$\wt{A}_{1}^{}$-configuration ($\wt{A}_{1}^{*}$-configuration) 
such that their intersection takes 
place on~$R$ with multiplicity~1 or~2 (below referred to as type~1 
and~2, respectively).

\begin{The} \label{thirdtheorem}
    $i)$ There exist $\wt{E}_{8}$-special $\mathbf{Z}/2$- 
and $\pmb{\alpha}_{2}$-surfaces. They are all classified by the 
preceding two theorems. 

$ii)$ There exist no $(\wt{A}_{1}^{}+\wt{E}_{7})$-special Enriques 
surfaces. All $(\wt{E}_{7}+R)$-surfaces are 
$(\wt{A}_{1}^{*}+\wt{E}_{7})$-special. They are all classified by the 
preceding two theorems. In particular,  
there exist $\pmb{\alpha}_{2}$-surfaces 
that are $(\wt{A}_{1}^{*}+\wt{E}_{7})$-special. 
They are all of type~2. And there exist 
$\mathbf{Z}/2$-surfaces that are 
$(\wt{A}_{1}^{*}+\wt{E}_{7})$-special of both types. 
\end{The}

\ned
I would like to thank Torsten Ekedahl for valuable help during the 
preparation of this paper. It was he who 
suggested to me the idea of studying 
Enriques surfaces using this particular kind of construction. 

\section{Proof of theorem~\ref{firsttheorem}}
Recall first the technique of resolving the singularities on a 
surface~$X'$ doubly covering another surface~$Y$ by resolving the 
singularities of the branch locus on~$Y$ and thereby to 
some extent being able to forget
about the third dimension that the double covering lives in. In 
characteristic two it is no longer appropriate to speak about the  
branch locus, but the basic principle is the same. For instance, to 
resolve (in characteristic two) $z^2=xy$ one would for instance put 
$x{\mapsto}x_{1}y$ to transform the RHS into $x_{1}y^2$, but 
then we must recall the third dimension and put $z{\mapsto}z_{1}y$ to 
obtain the equation $z^2_{1}y^2=x_{1}y^2$ and then removing the 
exceptional plane $\{y^2=0\}$ gives the new and smooth surface 
$z^2_{1}=x_{1}$ in local coordinates. In practise one never bothers 
about introducing the new symbol~$z_{1}$ but simply keeps~$z$ as third 
variable, dividing through by the square of the old variable~$y$ after 
each blowing up. It may very well happen that the blowing up 
introduces a \emph{curve} of singularities on the new surface. This 
corresponds exactly to the case when it is possible to divide the new 
RHS by \emph{higher} powers of the square~$y^2$ of the old variable. 
Performing this reduction amounts the normalisation of the surface, 
and the arithmetic genus drops by one for each extra 
power of~$y^2$ that it is possible to divide by. 
This and other facts about double-cover 
singularities can be found in~\cite{BPV}. 

In characteristic two the procedure takes on a somewhat unfamiliar 
appearance in the \emph{separable} case, that is, the case when the 
double covering is given by an equation $z^2+zg(x,y)+f(x,y)=0$ 
with~$g$ not identically zero. The 
differences are the following: $i)$~the singularities 
of the surface cannot be interpreted simply as the points that map to 
the singularities of the curve $f=0$; $ii)$~after the blowing up we 
shall divide the new version of~$f$ by~$y^2$ as before, but the new 
version of~$g$ should be divided by~$y$. 

We denote an unresolved surface as given by the equations 
by~$X'$. It is a double covering of a $\mathbf{P}^1$-bundle~$Y$  
over~$\mathbf{P}^1$ where the coordinates on the latter are~$x$ 
and~$y$ and the coordinates on a fiber of $Y{\ra}\mathbf{P}^1$ are~$s$ 
and~$t$. We let~$X''$ be the surface obtained by performing a minimal 
resolution of the singularities of~$X'$. It is not a minimal model, 
that is, it contains some smooth rational curves with 
self-intersection~$-1$. We denote the minimal model by~$X$. To 
sum up, we have the diagram
\begin{gather}\label{xypic.diagram}
\xymatrix@C+0pc@R+0pc{
  X' \ar[dr]^\phi  & X'' \ar[l]\ar[d]^\psi \ar[r] & X \ar[d]^\pi \\
                 & Y \ar[r]                  & \mathbf{P}^1
        } 
\end{gather} 
where~$\pi$ is a curve fibration with genus-one curves as fibers. 

There is a description of $\mathbf{P}^1$-bundles over $\mathbf{P}^1$ 
in terms of multihomogeneous coordinates. Indeed, rational functions 
on such a surface are given by quotients of polynomials in four variables 
satisfying two homogeneity conditions. In our case we have in the 
previous notation 
$d_{s}+d_{t}=\mathit{const.}$ and 
$2d_{t}+d_{x}+d_{y}=\mathit{const.}$, 
and the surface~$Y$ is thus seen to be the rational ruled surface with 
a section $C_{s}$ with self-intersection $-2$. This surface has,
cf.~\cite[V.2.10]{Har}, canonical class equal to $-2C_{t}$, 
writing~$C_{t}$ for the curve that is cut out 
by $t=0$, etcetera and its~$\chi$ equals~1 since it is rational. 

The surface~$X'$ has canonical class the pullback of that of~$Y$ plus 
that of the line bundle~$L^{-1}$ over~$Y$ that defines the double cover 
and as $L^{-1}\cong\mathcal{O}(2C_{t}+C_{x})$ the canonical class of~$X'$ 
is~$C_{x}$, that is, a fiber of $X'{\ra}\mathbf{P}^1$. 
Write~$D$ for the divisor $2C_{t}+C_{x}$. The~$\chi$ of the line 
bundle~$L$ is by the Riemann--Roch formula for surfaces equal to 
$\frac{(-D)^2-(-D)\cdot{K_{Y}}}{2}+1$, which equals 
$\frac{12-10}{2}+1=2$. This gives~$\chi(X')=3$, since it is 
the sum of that of~$L$ and that of~$Y$ 
by the discussion at the beginning of~\cite{CD}. 

There remains to prove the following. First that the surfaces have 
only isolated singularities. Next that the resolution 
$X''{\ra}X'$ subtracts two half-fibers from the canonical class and 
that the arithmetical genus drops from~3 to~1. This ensures that~$X$ 
is an Enriques surface. Thirdly we must of course prove that the 
surfaces have the stated curve configurations, more precisely that 
they arise as a resolution graph. 

As for isolated singularities, a curve of singularities must meet 
either the section~$R$ or a fiber of $X'{\ra}\mathbf{P}^1$, 
and one can readily see that all the given surfaces are 
nonsingular along~$R$ and along a typical fiber. 

As for the two invariants, the claim will follow if~$X'$ has exactly 
two genus one singularities or one of genus two, where the genus of a 
singularity is the amount by which the arithmetic genus of the surface 
drops during resolution. We know in advance that~$X'$ can't have more 
than that, since the~$\chi$ cannot be strictly less than~$1$, which 
is due to the fact that it is non-trivial as a genus-one fibration 
(that is, it does not split as a direct product), and then 
the first Betti number~$b_{1}$ vanishes, cf.~\cite[Cor.5.2.2]{CD}, 
which implies $\chi\geq1$. As for the canonical class 
the genus-one singularities will be seen to subtract a half-fiber each 
from the canonical class, so that it will add up to the difference 
between two distinct half-fibers, which is the classical situation on 
an Enriques surface, and in the other case, that 
of~$\pmb{\alpha}_{2}$ type surfaces, the genus two singularity occurs 
at the wild fiber, subtracting thus a fiber and making the canonical 
class trivial, as it should. 

I proceed to study the resolution process in the $\mathbf{Z}/2$ case. 
The two double fibers occur over the points $x=0$ and $y=0$ 
on~$\mathbf{P}^1$. The two nonrational singularities occur when we in 
addition put $t=0$. The curve $C_{s}=\{s=0\}$ is the wished-for 
curve~$R$. It is a double section, that is, $R{\ra}\mathbf{P}^1$ is a 
degree-two map. The two lines that make up the double fibers can be 
contracted after the resolution, that is, they will have 
self-intersection $-1$. Doing this will connect~$R$ to the resolution 
graph so that we obtain the stated configuration. 

The nonrational singularity at $y=0$ is of no further importance for 
us. For the one at $x=0$ the following happens. The first blowing up 
produces on the modified~$X'$ a rational 
curve as exceptional divisor. It 
has an isolated nonrational singularity on it, and 
resolving this we would have obtained an 
elliptic curve as exceptional divisor if we had used a somewhat more 
generic choice of parameters than in the theorem. In any case the 
original fiber over $x=0$ can be blown down as well as the first 
obtained exceptional divisor. The second exceptional divisor has, 
using the stated parameters, a rational singularity on it, and one 
can check that it is of the appropriate non-extended Dynkin type by 
continuing the resolution process. 

In the~$\pmb{\alpha}_{2}$ case there is a genus two singularity at 
$x=t=0$. In the most special case we have locally 
upon putting $a_{80}=1$ the equation $z^2+x^3+x^8t+xt^4=0$. 
Blowing up by $x{\mapsto}yt$ gives after normalisation 
$y^3t+y^8t^7+yt^3$ with singularity at the origin. 
Blowing up by $y{\mapsto}wt$ gives after normalisation 
$w^3+w^8t^{11}+w$ with singularity at $w+1=t=0$. This time 
the normalisation made the arithmetic genus of the double cover drop. 
Changing $w{\mapsto}\ol{w}+1$ gives after removal of squares 
$\ol{w}^3+t^{11}+\ol{w}^8t^{11}$ which one can approximate by 
$\ol{w}^3+t^{11}$ without influencing the resolution process, as one 
easily verifies. Blowing up by $\ol{w}{\mapsto}ut$ gives after 
normalisation $u^3t+t^9$ with singularity at the origin. 
Blowing up by $u{\mapsto}st$ gives after 
normalisation $s^3+t^5$. This time the arithmetical genus drops for 
the second time. We may thus expect a genus-zero singularity this 
time, and more precisely we have $z^2+s^3+t^5$, which is exactly the 
standard formula for the~$E_{8}$ surface singularity, cf.~\cite{CD} 
or~\cite{Lip}. The pieces of the resolution graph that appear before 
the curve that came into being when the genus dropped for the second 
time can be blown down, so that curve is the component in the 
extended Dynkin graph that will meet~$R$. The other cases are treated 
similarly. 

\section{Proof of theorem~\ref{secondtheorem}}
As mentioned the stated equations define surfaces~$X'$ which make up a 
curve fibration over~$\mathbf{P}^1$ with genus-one curves as fibers, 
and the mapping factors through a double covering of a 
$\mathbf{P}^1$-bundle~$Y$, so that the restriction to a typical fiber 
is an elliptic curve doubly folded over a line (except of course 
if the curve fibration is quasi-elliptic). 

\begin{Pro} An Enriques surface with an 
$(\wt{E}_{i}+R)$-configuration can be represented as a curve fibration 
over~$\mathbf{P}^1$ with the~$\wt{E}_{i}$ configuration as half-fiber 
and~$R$ as double-section. The relative projective mapping 
associated to the line bundle $\mathcal{O}_{X'}(R)$ is 
(after blowing up of base points) a degree-two mapping
to a~$\mathbf{P}^1$-bundle~$Y$ with the same twisting as the one 
in the previous section, and the double covering can be written by an 
equation of the same multihomogeneity as the ones in 
theorem~\ref{firsttheorem}.  
\end{Pro}
\bevis The existence of the genus-one fibering is Proposition 3.1.2 
(page 171) of the book~\cite{CD} by Cossec and Dolgachev. Then 
Theorem 4.4.1 (page 240) in \emph{loco citato} shows that the number 
of base points of the linear system is the same for all 
$(\wt{E}_{i}+R)$-surfaces. The calculation of the twisting of~$Y$ and 
the form of the equation depends only on this fact and on 
numerical invariants that are the same for all 
$(\wt{E}_{i}+R)$-surfaces. 

Indeed, for the twisting, after two blowing ups the self intersection 
of~$R$ is $-4$ and hence the exceptional section of~$Y$ must have 
self intersection~$-2$. 

The push-forward of the structure 
sheaf of a double covering surjects onto a line-bundle~$L$ with a 
trivial line bundle as kernel. The covering can be written in the 
stated form if this extension splits. We need to know that 
$H^1(L^{-1})=0$. 

The line bundle~$L$ is determined by the fact that the canonical 
bundle of the cover (which is known) is the pullback of 
$\omega_{Y}\otimes L$, cf.~\cite{CD} page~12. In our case this gives 
that $L^{-1}\cong\mathcal{O}(2C_{t}+C_{x})$ with the previous 
notation. The fact that $H^1(Y,L^{-1})=0$ follows from the Leray 
spectral sequence associated to $Y\ra\mathbf{P}^1$. 
\qed 

\ned 
The diagram~(\ref{xypic.diagram}) now takes on the 
meaning that the mapping $X''{\ra}X$ is a blowing up of base points  
of the relative linear system. The part of the fiber which is 
not the exceptional divisor is contracted by the map $X{\ra}X'$ 
giving rise to non-rational singularities and we know the location 
of these on~$X'$ by the explicit checking that we have done. 
There remains for us to prove that only the 
\emph{stated} equations are those that arise 
from Enriques surfaces of the given type, and also to check the 
uniqueness claim in the $\mathbf{Z}/2$ case. 

We have thus the general equation of the given multihomogeneity 
\begin{gather} \label{general.eqn}
z^2 + z(B_5s^2 + B_3st + B_1t^2) + \notag\\
A_{10}s^4 + A_8s^3t+ A_6s^2t^2 +  A_4st^3 + A_2t^4=0,  
\end{gather}
where the $A_i$'s and $B_i$'s are forms in $x,y$ of degree~$i$. 
Writing more compactly $z^2+zg+f$, where~$f$ and~$g$ are polynomials 
in the four variables, we first see that any such surface is 
isomorphic to one where~$f$ is square free. 
Namely, put $z\mapsto z+h$, where~$h$ is a polynomial of the same 
degree as~$g$ such that the square parts of 
$h^2+gh$ and~$f$ are equal. The following lemma shows that  
such a polynomial can be found. This reduction is 
of course entirely trivial if the characteristic isn't two. 

\begin{Lem}
Let~$u$ and~$g$ be polynomials in one or several variables 
over a field of characteristic~$p$ such that the multidegree of~$g$ 
is $p-1$ times that of~$u$. Then there is a (non-unique) 
polynomial~$h$ of the same degree as~$u$ 
such that the $p$'th power part of $h^p+gh$ is~$u^p$. 
\end{Lem}    
\bevis 
Let~$R$ be the space of polynomials of the same degree as~$u$. 
The proof is non-constructive. 
We will need to consider~$R$ as
an algebraic group, rather than as a 
vector space. Define a group scheme
endomorphism~$\tau$ by first sending 
an $h\in{R}$ to the $p$'th power part of
$h^p+gh$ and then dividing all exponents 
in the resulting polynomial by $p$. 
For a given~$g$ we have an endomorphism 
of~$R$. If we  know that the image contains any given~$u$, then a 
polynomial~$h$ in its preimage will satisfy the stated condition. 
We shall thus prove that~$\tau$ is surjective.  
Suppose $\mbox{ker}(\tau)$ is \emph{not} 
a finite group scheme. 
Then there exists \cite[sect. 20]{Hum} a non-constant
homomorphism $\sigma\colon\mathbf{G}_{a}\ra{R}$, 
such that composing with~$\tau$
kills it. Explicitly~$\sigma$ is given by 
a vector of additive polynomials:
\begin{equation*}
x\mapsto(a_{10}x+a_{11}x^p+a_{12}x^{p^2}+\cdots,
\,a_{20}x+\cdots,\,\dots)  ,.    
\end{equation*}
But since~$\tau$ is the sum of the Frobenius 
map and a linear transformation,
looking at the terms of highest degree that 
occur in the expression for~$\sigma$, it is impossible that the 
composition $\tau\sigma$ is the zero map.
So $\mbox{ker}(\tau)$ is finite, showing 
that~$\tau$ is surjective. 
\qed

\ned
We thus get rid of the monomials in~$A_{10}$,~$A_{6}$ and~$A_{4}$ 
that have even degree in~$x$ or~$y$. From now on the two 
types~$\mathbf{Z}/2$ and~$\pmb{\alpha}_{2}$ must be treated 
separately. We begin with the former. 

\subsection{The $\mathbf{Z}/2$ case}
We may put the two non-rational singularities that arise from the 
base-points at $x=t=0$ and $y=t=0$. To get double fibers~$x$ and~$y$ 
must divide~$f$ and~$g$. The equation now looks like 
\begin{gather*}
z^2 + z\bigl((B_{14}+B_{23}+B_{32}+B_{41})s^2+
 (B_{12}+B_{21})st\bigr)  \\
 (A_{19}+A_{37}+A_{55}+A_{73}+A_{91})s^4+\\
 (A_{17}+A_{26}+A_{35}+ A_{44}+A_{53}+A_{62}+A_{71})s^3t+\\
 (A_{15}+A_{33}+A_{51})s^2t^2+
 (A_{13}+ A_{22}+A_{31})st^3+A_{11}t^4=0.  
\end{gather*}
where $A_{ij}$ denotes $a_{ij}x^{i}y^{j}$, where 
$a_{ij}$ are parameters, and similarly for $B_{ij}$. 

To have a singularity at $x=t=0$ we must have $A_{19}=0$. To get
a non-rational singularity, the quadratic part of the
polynomial must be a square according to~\cite{Lip}. 
Therefore $B_{14}=A_{17}=0$. The cubic 
part of~$f$ must be a cube, cf.~\cite{Lip}, giving 
$A_{26}=A_{15}=0$. To have~$X'$ smooth over the rest of~$C_x$, which 
it must be since this is the exceptional divisor of the blowing-up of 
a base-point, we must have $A_{13}=0$ and $A_{11}\neq0$, so 
put $A_{11}=1$ by multiplying the $z$-variable by a 
constant. Moreover~$A_{37}$ must be nonvanishing 
since otherwise the following 
happens: the genus drops at the first blowing up and then there is at 
least one more blowing up to do along the old fiber, which gives 
it the wrong self-intersection. 
So we may put $A_{37}=1$ by using the third
automorphism of $\mathbf{P}^1$ or $t\mapsto\lambda t$. 
After repeating everything at the other point,
we get the desired polynomial by putting $A_{55}=0$, which is 
possible by using
$t\mapsto t+\lambda sxy$. Doing this, a square is re-introduced 
into~$f$, namely a linear combination of $x^4y^6s^4$ and 
$x^6y^4s^4$. This can however be removed by 
$z\mapsto z+\mu_1 x^2y^3s^2 +\mu_2 x^3y^2s^2$. 
The equation now looks like 
\begin{gather}\label{any.Z2}
z^2 + z\bigl((B_{23}+B_{32})s^2+ (B_{12}+B_{21})st\bigr)  \\
 (x^3y^7+x^7y^3)s^4+(A_{35}+ A_{44}+A_{53})s^3t+
 A_{33}s^2t^2+A_{22}st^3+xyt^4=0.  \notag
\end{gather}
and we can make the uniqueness claim that two different equations of 
this type give non-isomorphic surfaces since we have used up all 
automorphisms of the ruled surface~$Y$. 

\ned
We proceed to determine conditions on the parameters in order to get 
a~$\wt{E}_{6}$ singularity at $x=t=0$. We resolve the singularity 
there, and for each blowing up we get some conditions in order that 
the resolution process continues. We shall see that we successively 
get the following conditions: 
\begin{align}
    b_{12}=a_{35}+a_{22}&=0 \label{cond.i}\\
    b_{23}=a_{22}&=0         \label{cond.ii}\\
    a_{44}+\sqrt{a_{33}}b_{22}=b_{22}+a_{33}&=0\label{cond.iii}
\end{align}
which together imply the equation stated for a $\mathbf{Z}/2$ type 
with $\wt{E}_{6}+R$ configuration in theorem~\ref{firsttheorem} if we 
introduce the auxiliary parameter~$v$, where 
$v^2=b_{22}=a_{33}$ and $v^3=a_{44}$.  

Blow up once. The genus does not drop and there is still a 
singularity at the old fiber. The next blowing up makes the 
genus drop and with a generic choice of parameters we get a smooth 
elliptic curve as exceptional divisor on the double cover. We must 
have a further singularity and we get that upon 
imposing~(\ref{cond.i}). 
The singularity is of type~$D_{n}$ or~$E_{n}$ if we also 
impose~(\ref{cond.ii}). It is of~$E_{n}$ type if we also 
impose~(\ref{cond.iii}). To check 
these statements one may use Lipman's 
conditions~\cite{Lip}.

Thus far we have the equation for a $\mathbf{Z}/2$ type surface 
with $\wt{E}_{6}+R$ configuration. To obtain a $\wt{E}_{7}+R$ 
configuration instead one must put $v=0$. 
To obtain a $\wt{E}_{8}+R$ 
configuration instead one must put~$b_{32}=0$ and then we must 
require~$a_{53}\neq0$ to avoid a non-isolated singularity, but if 
this condition is satisfied we indeed get~$\wt{E}_{8}$. 

\subsection{The $\pmb{\alpha}_{2}$ case}
This time there is only one non-rational singularity, but is has 
\emph{genus two}. Let us put it at $x=t=0$. 
To get a double fiber at $x=0$, we see that~$x$ 
must divide~$f$ and~$g$. We get conditions on the parameters 
in order that there be a genus two singularity at the 
prescribed spot. 

The equation that is analogous to~(\ref{any.Z2}) is in 
the~$\pmb{\alpha}_{2}$ case the following one: 
\begin{gather*}
z^2 + z\bigl((w'x^3y^2+B_{41}+B_{50})s^2+B_{30}st\bigr)  \\
 x^3y^7s^4+ (w'x^4y^4+wx^5y^3+A_{62}+A_{71}+A_{80})s^3t+ \\
 A_{51}s^2t^2+(w'x^3y+wx^4)st^3+xyt^4=0\,,  
\end{gather*}
where~$w$ and~$w'$ are parameters just like the~$a_{ij}$'s 
and~$b_{ij}$'s. These equations thus describe 
Enriques surfaces of $\pmb{\alpha}_{2}$ 
type without any condition that they should have some particular 
configuration of rational curves on them. The proof runs in parallel 
to the $\mathbf{Z}/2$ case. I shall not present the details, but at 
least I shall describe how the resolution process runs. The genus 
drops by one at the second blowing up just as in the 
$\mathbf{Z}/2$ case and the first exceptional 
divisor can be blown down. The singularity on the second 
exceptional divisor is this time a genus one singularity and the 
genus drops at the next blowing up. If the parameters are general,  
the new curve is smooth. To obtain a singularity of $E_{6}$ type 
on it we shall impose 
\begin{equation*}
    w'=b_{30}=b_{41}+vw=a_{62}=a_{71}=a_{51}=0 
\end{equation*}
which gives 
\begin{gather*}
z^2 + z(wx^4y+B_{50})s^2 +
 x^3y^7s^4+ (wx^5y^3+A_{80})s^3t+ wx^4st^3+xyt^4=0   
\end{gather*}
which is the equation stated in theorem~\ref{firsttheorem} for the 
$\pmb{\alpha}_{2}$ type with an $\wt{E}_{6}+R$ configuration. Putting 
$w=0$ it degenerates into a~$\wt{E}_{7}$ and putting in addition 
$b_{50}=0$ it degenerates into a~$\wt{E}_{8}$, except if $a_{80}=0$, 
in which case the singularity ceases to be isolated. 

\subsubsection{Putting a parameter equal to one in the 
$\pmb{\alpha}_{2}$ case}
We have claimed that one can put~$w$,~$b_{50}$ or~$a_{80}$ 
equal to~1 in the $\wt{E}_6$, $\wt{E}_7$ and $\wt{E}_8$ cases, 
respectively. One uses variable changes $x\mapsto \lambda x$ and 
similarly for $y$, $s$, $t$ and~$z$. The isomorphism type of the 
surface remains the same. And we can also rescale the entire equation 
by a multiplicative constant~$k$. 

Let us consider the~$\wt{E}_6$ case, which is most difficult. 
There are three monomials whose coefficient is required to be~1. They 
give rise to the conditions $k\lambda_{z}^2=1$, 
$k\lambda_{x}^3\lambda_{y}^7\lambda_{s}^4=1$ and 
$k\lambda_{x}\lambda_{y}\lambda_{t}^4=1$. Here we have 
written~$\lambda_{z}$ for the rescaling constant of~$z$, etcetera. 
There are three monomials whose coefficient is~$w$. They 
give rise to the conditions 
$k\lambda_{z}\lambda_{x}^4\lambda_{y}\lambda_{s}^2=w^{-1}$, 
$k\lambda_{x}^5\lambda_{y}^3\lambda_{s}^3\lambda_{t}=w^{-1}$ and
$k\lambda_{x}^4\lambda_{s}\lambda_{t}^3=w^{-1}$. We have thus 
six equations and six unknown variables. We can use the evident 
multiplicative version of Gauss elimination. A solution is given by 
$\lambda_{y}=\lambda_{s}=1$, $\lambda_x=w^{-2/5}$, 
$\lambda_t=w^{-1/5}$, $\lambda_z=w^{-3/5}$ and $k=w^{6/5}$. 

In the~$\wt{E}_7$ case we get in similar fashion the conditions
$k\lambda_{z}^2=1$, $k\lambda_{x}^3\lambda_{y}^7\lambda_{s}^4=1$, 
$k\lambda_{z}\lambda_{x}^5\lambda_{s}^2=b_{50}^{-1}$, 
$k\lambda_{x}\lambda_{y}\lambda_{t}^4=1$ and a solution
$k=\lambda_{z}=\lambda_{x}=1$, $\lambda_{s}=b_{50}^{-2}$, 
$\lambda_{y}=b_{50}^{-2/7}$ and $\lambda_{t}=b_{50}^{-1/14}$. 

In the~$\wt{E}_8$ case we get the conditions 
$k\lambda_{z}=1$, $k\lambda_{x}^3\lambda_{y}^7\lambda_{s}^4=1$, 
$k\lambda_{x}^8\lambda_{s}^3\lambda_{t}=a_{80}^{-1}$ and a solution
$k=a_{80}^{1/11}$, $\lambda_{z}=a_{80}^{-1/22}$, 
$\lambda_{s}=\lambda_{t}=1$, $\lambda_{x}=a_{80}^{-3/22}$ and 
$\lambda_{y}=a_{80}^{1/22}$.

\section{Proof of theorem~\ref{thirdtheorem}}
The $\wt{E}_{8}$-special case is contained in 
theorem~\ref{firsttheorem}. 

The fibering associated to an 
$\wt{E}_{7}+R$ or $\wt{E}_{8}+R$ configuration is quasi-elliptic. 
This can be read off directly from the equations. Please notice that 
quasi-ellipticity does not imply that the coefficients in front of~$z$ 
vanish. 

Due to quasi-ellipticity, it follows from the Euler characteristic 
formula for a curve pencil (\cite{CD} page 290) that any 
$(\wt{E}_{7}+R)$-surface is $(\wt{A}_{1}^{*}+\wt{E}_{7})$-special, but 
in the $\mathbf{Z}/2$ case it is not clear if it is of 
type~1 or~2. One can however check 
explicitly that both types can occur. The formula says that the 
($\ell$-adic) Euler characteristic is the product of those of the base 
and a typical fiber plus contributions from degenerate fibers. In the 
quasi-elliptic case the local contributions are nothing more than the 
amount that the Euler characteristic of the fiber exceeds the typical 
value. The Euler characteristic of an Enriques surface is~12. The 
base has Euler characteristic 2, and the same holds for a typical 
fiber in the quasi-elliptic case, so the local contributions must sum 
up to~8. The~$\wt{E}_{7}$ fiber gives~7, and hence there must be 
exactly one more extra fiber component, and that fiber must then be 
of~$A_{1}^{*}$ type, that is, two smooth rational curves intersecting 
with multiplicity two.


\end{document}